\documentclass[12pt]{article}
\usepackage{amsmath,amssymb,amsfonts}

\newcommand{\lo}{\longrightarrow}
\newcommand{\B}{\beta}
\newcommand{\var}{\varepsilon}
\newcommand{\st}{\stackrel}
\newcommand{\pa}{\paragraph}
\begin{document}
\begin{center}
{\bf The Dynamics of a Vertically Transmitted Disease }\\[1cm]
M.R. Razvan \\

{\it \small Institute for Studies in Theoretical Physics and Mathematics\\
\small P.O.Box: $19395-5746$, Tehran, IRAN} \\
Email: razvan@karun.ipm.ac.ir \\
Fax: 009821-2290648
\end{center}

\vspace*{5mm}
\begin{abstract}
An SIRS epidemiological model for a vertically transmitted disease
is discussed. We give a complete global analysis in terms of
three explicit threshold parameters which respectively govern the
existence and stability of an endemic proportion equilibrium, the
increase of the total population and the growth of the infective
population. This paper generalizes the results of Busenberg and
van den Driessche.
\end{abstract}


\section{Introduction}
In 1990, a complete global analysis of an SIRS model of disease
transmission in a population with varying size was given by
Busenberg and van den Driessch \cite{Bv1}. In an SIRS
epidemiological model, we divide the population into three
groups, Susceptibles, Infectives and Removeds, and the problem is
to examine the behavior of the size of each group when the time
goes to infinity. They considered a disease with horizonal
transmission, that is a disease which is transmitted by contact
between an infective and a susceptible individual. We intend to
generalize their results for a vertically transmitted disease,
that is a disease which is also transmitted from infective parents
to their newborns. The assumption of vertical transmission has two
consequences. First it causes some newborns to die and forces us
to assume that the birth rate differs from one group to another.
(See the demographic assumption in \cite{BH,MAM1,MAM2}.) The
second fact is that some newborns are infected, hence a group of
newborns enter to the infective class \cite{BC,BCP,C}. We also
assume that a part of these infected newborns are known and
removed after their birth. We shall show that the latter
parameter can play an important role in the epidemic process.

We consider an SIRS epidemiological model for a vertically
transmitted disease. In our model, the incidence function is of
proportionate mixing type introduced by Nold \cite{N}. Natural
births and deaths are assumed to be proportional to the class
numbers with different rates. We also assume that a proportion of
the infected newborns are known and removed. We follow \cite{Bv1}
to examine our model equation which is homogeneous of degree one.
We consider the proportions system and show that this system has
no periodic orbit in its feasibility region. This reduces our
analysis to the discussion of existence and stability of rest
points of a palnar system. The technique used here to show the
uniquness of endemic equilibria is based on the Poincar\'e index.
This technique has no hard analysis and can be easily applied to
other similar systems \cite{R,RK}. The reader can verify that our
results hold for similar SIRI systems as well \cite{Dv}.

We first in the next section state the model and a result
concerning the non-existence of certain types of solutions
\cite{Bv2}. We consider the proportions system and prove that
every solution in the feasibility region tends to a rest point of
this system. In Section 3., we introduce a quadratic planar
system with the same dynamics as the proportions system and then
we discuss the existence and stability of rest points of this
quadratic planar system. This gives a complete global analysis of
the proportions system which is used to provide a global analysis
of the original system is Section 4.

\section{The model}
In order to derive our model equations, we divide the population
into three classes, the susceptible, the infective and the
removed individuals with total numbers $S$, $I$ and $R$
respectively,. We set $N=S+I+R$ which is the total size of the
population. The following parameters are used in our model
equations:
\begin{itemize}
\item[$b_0$:] per capita birth rate of susceptible individuals,
\item[$b_1$:] per capita birth rate of infective individuals born uninfected,
\item[$\B$:]   per capita birth rate of infective individuals born infected,
\item[$b_2$:] per capita birth rate of removed individuals,
\item[$d$:]   per capita disease free death rate,
\item[$\var_1$:] excess per capita death rate of infected individuals,
\item[$\var_2$:] excess per capita death rate of removed individuals,
\item[$\alpha$:] per capita removal rate of infective individuals,
\item[$\gamma$:] per capita recovery rate of removed individuals,
\item[$\lambda$:] effective per capita contact rate of infective individuals.
\end{itemize}
As mentioned before, we assume that the infected newborns enter
the classes $I$ and $R$ of proportions $\B_1$ and $\B_2$
respectively, hence $\B=\B_1+\B_2$. In this paper, all the above
parameters, are positive, however some of them can also be zero.
These hypotheses yield the following system of differential
equation in $\mathbb{R}_3^+$, where `` $'$ '' denotes the
derivatives with respect to $t$, the time.

  $\left\{
\begin{array}{l}
{S'}=(b_0-d)S+b_1I+(b_2+\gamma)R-\lambda
\dfrac{IS}{N}\hspace*{6cm}
 \hfill{(2-1)}\\
{I'}=(\B_1-d-\var_1-\alpha)I+ \lambda
\dfrac{IS}{N} \hfill{(2-2)}\\
{R'}=(\B_2+\alpha)I-(d+\var_2+\gamma)R
\hfill{(2-3)}\end{array}\right.$

\noindent where $\lambda \dfrac{IS}{N}$ is of the proportionate or
random mixing type \cite{Hv,N}. By adding the above three
equations, the total population equation is
$$N=b_0S+(b_1+\B-\var_1)I+(b_2-\var_2)R-dN.$$
If we consider the proportions $s=\dfrac{S}{N}$, $i=\dfrac{I}{N}$,
$r=\dfrac{R}{N}$, we get the following system of equations.

$\left\{ \begin{array}{l} s'=b_0 s+b_1 i+(b_2+\gamma)r-b_0
s^2-(b_1+\B+\lambda-\var_1)is
-(b_2-\var_2) sr\hspace*{1cm} \hfill{(2-1)'}\\
i'=(\B_1-\var_1-\alpha)i+(\lambda-b_0) is -(b_1+\B-\var_1) \
i^2 -(b_2-\var_2) \ ir \hfill{(2-2)'}\\
r'=(\B_2+\alpha)i-(\var_2+\gamma)r-b_0 sr-(b_1+\B-\var_1) ir
-(b_2-\var_2)r^2 \hfill{(2-3)'}
\end{array} \right.$

\noindent This system is called proportions system and the
feasibility region of this system  is the triangle
$$D=\{(s,i,r)|s\geq 0, i\geq 0, r\geq 0, s+i+r=1\}.$$
If we set $\Sigma=s+i+r$ then
$\Sigma'=(1-\Sigma)(b_0s+(b_1+\B-\var_1)i+(b_2-\var_2)r)$. Thus
the plane $\Sigma=1$ is invariant. Moreover, on the sides of $D$,
we have:

\begin{center}
 $\left\{ \begin{array}{l}
s=0 \Rightarrow s'= b_1 i+(b_2+\gamma)r \geq 0,\\
i=0 \Rightarrow i'=0 \quad  \mbox{hence the line}\  \{i=0\} \  \mbox{is invariant}, \\
r=0 \Rightarrow r'=(\B_2+\alpha)i.
\end{array}\right .$
\end{center}

\noindent Therefore, $D$ is positively invariant. On the invariant
line $\{\Sigma=1\}\cap \{i=0\}$, we have
$$r'=-(\var_2+\gamma)r-b_0
r(1-r)-(b_2-\var_2)r^2=-(b_0+\var_2+\gamma)r+(b_0 -b_2+\var_2)
r^2.$$ It follows that this invariant line contains two rest
points, the Disease-Free Equilibrium $(1,0,0)$ and possibly
 another one which is outside of D. It is easy to see
that the DFE attracts the side $D\cap \{i=0\}$. Furthermore, our vector field
is strictly inward on the other sides of $D$. Thus $\st{\circ}{D}$,
 the interior of $D$, is positively invariant too.
The following theorem reduces our problem to the discussion of
existence and stability of rest points in $D$.

\pa{Theorem 2.1.} The $\omega$-limit set of any solution for the
system $(2-1)'-(2-3)'$ with initial point in $D$ is a rest point
in $D$.

\pa{Proof.} Since $D$ is compact and positively invariant, the
$\omega$-limit set of any solution with initial point in $D$ is a
compact nonempty invariant subset of $D$. Here we use the
Poincar\'e-Bendixon theorem and the terminology used to prove it
 \cite{P}. In the next section, we will see that this system has at most
three rest points in $D$. Thus it satisfies the assumptions of
Poincar\'e-Bendixon theorem. We follow \cite{Bv1,Bv2} and define
the vector field $g=(g_1,g_2,g_3)$ on $D$ by
\begin{align*}
g_1(i, r)=\left [0, -\dfrac{f_3(i, r)}{ir}, \dfrac{f_2(i, r)}{ir}\right ] , &\\
g_2(s, r)=\left [\dfrac{f_3(s,r)}{sr}, 0, -\dfrac{f_1(s, r)}{sr}\right ], &\\
g_3(s, i)=\left [-\dfrac{f_2(s, i)}{si}, \dfrac{f_1(s, i)}{si},
0\right ], &
\end{align*}
where $f_1,\ f_2$ and $f_3$ are the right hand side of $(2-1)'$,
$(2-2)'$ and $(2-3)'$ reduced to functions of two variables by
using $\sum =1$ respectively.
 Clearly $g.f=0$ in $\st{\circ}{D}$ and
after some computations \cite{Dv}, we get
$$(\mbox{curl}\ g).
(1,1,1)=-\left(\frac{b_1}{s^2r}+\frac{b_2+\gamma}{s^2i}+\frac{\B_2+\alpha}
{sr^2} \right).$$ Since the DFE is the only invariant subset of
$\partial
 D$ (i.e. the boundary of
$D$), the $\omega$-limit set must have some regular point in
$\st{\circ}{D}$
 if it is
not a rest point. Let $x$ be such a regular point and $h$ be the
first return map (Poincar\'e map) defined on a tranversal at $x$.
For a point $y$ near $x$ on the transversal, Let $V$ be the region
surrounded by the orbit $\Gamma$ from $y$ to $h(y)$ and the
segment between them. (This region is known as Bendixon Sack, Sec
Fig 2.1) Now by Stokes' theorem
$$\int\int_V (\mbox{curl}\
g).(1,1,1) d\sigma=\int_{\Gamma} g.f dt+\int_0^1 g(ty+(1-t)h(y)).
(y-h(y))dt$$
Since $g.f=0$ and $h(x)=x$, the right hand side of
the above equality tends to zero when $y$ tends to $x$, but the
left hand side tends to the integral over the region bounded by
the $\omega$-limit set, this is a contradiction since
$(\mbox{curl}\ g).(1,1,1)<0$ in $\st{\circ}{D}$. $\square$
\begin{center}
{\large {
\unitlength 1.00mm
\linethickness{1.1pt}
\begin{picture}(77.33,80.67)
\put(35.33,80.67){\line(1,-1){17.33}}
\bezier{196}(68.67,63.67)(77.33,44.33)(50.00,38.00)
\bezier{100}(28.33,57.33)(27.67,67.67)(37.67,78.00)
\bezier{156}(50.00,38.00)(31.00,37.33)(28.33,56.67)
\put(35.00,75.00){\vector(1,1){1.00}}
\bezier{136}(68.67,63.67)(59.67,79.33)(46.67,69.67)
\put(49.33,71.33){\vector(3,2){1.00}}
\put(45.00,30.00){\makebox(0,0)[cc]{{\small Fig. 2.1. The Bendixon
Sack.}}}
\end{picture}}}
\end{center}
\vspace*{-3cm}
\section{The planar system}
using the relation $s+i+r=1$, we see that our system is essentially two
dimensional. Thus we can eliminate one of the variable to arrive
at the following quadratic planar system:\\
$$\left\{
\begin{array}{l}
i'=(\lambda+\B_1-b_0-\var_1-\alpha)i+ (b_0+\var_1-\lambda
-b_1-\B)i^2+(b_0+\var_2-\lambda-b_2)ir\quad  \hfill{(3-1)}\\
r'=(\B_2+\alpha)i-(b_0+\var_2+\gamma)r+(b_0+\var_1-b_1-\B)ir+(
b_0+\var_2-b_2)r^2 \hfill{(3-2)}
\end{array}\right.$$
The dynamics of the system $(2-1)'-(2-3')$ in $D$ is equivalent to the
dynamics of this planar system in the positively invariant region
$D_1=\{(i,r): i\geq 0,r\geq 0,i+r\leq 1\}$.
This quadratic system has at most four rest points and since $D_1$
misses a rest point  on $i=0$,
 there are at most three rest points in $D_1$. One of these
rest points is the origin which comes from the DFE. The matrix of the
linearization of the system $(3-1),(3-2)$ at the origin is:
$$\left[ \begin{array}{lc}
\lambda+\B_1-b_0-\var_1 -\alpha & 0\\
\B_2+\alpha & -(b_0+\var_2+\gamma)
\end{array}\right]$$
with the eigenvalues $\lambda+\B_1-b_0-\var_1-\alpha$ and
$-(b_0+\var_2+\gamma)$. Now we define the first threshold
parameter $R_0=\dfrac{\lambda+\B_1}{b_0+\var_1+\alpha}$ which
governs the stability of the origin.

\pa{Theorem 3.1.} The origin is globally asymptotically stable in
the feasibility region $D_1$ when $R_0\leq 1$ and it is a saddle
point when $R_0>1$.

\pa{Proof.} By Theorem 2.1. it is enough to prove that if $R_0\leq
1$, then the origin is the only rest point in $D_1$. If there
exists a rest point in $\st{\circ}{D_1}$, we have $i'=0$ and
$i\neq 0$ at this point. Thus it belongs to the line
$$(\lambda+\B_1-b_0-\var_1-\alpha)+(b_0+\var_1-\lambda-b_1-\B)i+
(b_0+\var_2-\lambda-b_2)r=0. \eqno{(3-3)} $$
 Since $s'=0$ at this point, from $(2-1)'$ we obtain
$$b_0 s+b_1 i+(b_2+\gamma)r -b_0 s^2-(b_1+\B+\lambda-\var_1) is
-(b_2-\var_2)sr=0.$$ and by using the relation $s+i+r=1$, we can
write
$$b_1i+(b_2+\gamma)r+(b_0+\var_1-b_1-\B-\lambda) is
+(b_0+\var_2-b_2)sr=0.$$ Multiplying (3-3) by $(-s)$ and adding
it to the above expression, we get the following equality
$$b_1i+(b_2+\gamma)r+\lambda sr+ (b_0+\var_1+\alpha-\lambda-\B_1)s=0.$$
But the left hand side is positive when
 $R_0\leq 1$ and this is a contradiction.
$\square$

When $R_0>1$, the origin is a saddle point and it does not attract
any point of $D_1-\{i=0\}$. Thus the orbits with initial point in
$D_1-\{i=0\}$ must be attracted by some rest points in
$\st{\circ}{D}_1$ by Theorem 2.1. These rest points belong to the
line (3-3) and the conic section
$$r'=(\B_2+\alpha)i-(b_0+\var_2+\gamma)r+(b_0+\var_1-b_1-\B)ir+
(b_0+\var_2-b_2)r^2=0.$$ It follows that there are at most two
rest points in $\st{\circ}{D}_1$. Notice that a nondegenerate
rest point of the planar system is obtained by a transverse
intersection of the line (3-3) and the above conic section. The
following lemma has two immediate consequences which will be very
helpful.

\pa{Lemma 3.2.} The trace of the linearization of the system
(3-1),(3-2) at a rest point in $\st{\circ}{D}_1$ is negative.

\pa{Proof:} We compute the trace at a rest point.
$$\begin{array}{ll}
\dfrac{\partial i'}{\partial i}& =(\lambda
+\B_1-b_0-\var_1-\alpha)+2(b_0+\var_1-\lambda
-b_1-\B)i+(b_0+\var_2 -\lambda-b_2) r,\\  \dfrac{\partial
r'}{\partial r} &=-( b_0+\var_2+\gamma)+(b_0+\var_1-b_1-\B)i+
2(b_0+\var_2-b_2)r.\end{array}$$ From $i'=0$ and $r'=0$ at a rest
point, we can write  $$\frac{\partial i'}{\partial
i}=(b_0+\var_1-\lambda-b_1-\B)i \quad \mbox{and} \quad
\frac{\partial
 r'}{\partial r}=(b_0+\var_2-b_2)r-(\B_2+\alpha)\frac{i}{r
}.$$ Since $s'=b_1i+(b_2+\gamma) r+(b_0+\var_1-\lambda -b_1-\B)
si+(b_0+\var_2-b_2) sr$, it follows that $(b_0+\var_1-\lambda
-b_1-\B)si+(b_0+\var_1-\lambda-b_1- \B) sr<0$ and hence
$\dfrac{\partial i'}{\partial i}+\dfrac{\partial r'}{\partial
r}<0$. $\square$

\pa{Corollary 3.3.} The system (3-1),(3-2) has no source point in
$\st{\circ}{D}_1$

\pa{Corollary 3.4.}  Every nondegenerate rest point in
$\st{\circ}{D}_1$ is hyperbolic.

\pa{Theorem 3.5.} If $R_0>1$, then there exists a unique rest point
$(i^*,r^*)$ in $\st{\circ}{D}_1$ which is hyperbolic and attracts
$D_1-\{i=0\}$.

\pa{Proof.} When $R_0>1$, the origin is a saddle point with the unstable
eigenvector
$$\left[\begin{array}{c}
\lambda+\B_1-b_0-\var_1-\alpha+b_0+\var_2+\gamma\\
\B_2+\alpha
\end{array}\right].$$
Since $R_0>1$, we have $\lambda+\B_1-b_0-\var_1-\alpha>0$ and
hence this vector belongs to the first quadrant of the plane
$(i,r)$. Since  $\st{\circ}{D}_1$ is positively invariant it
follows that a branch of the unstable manifold of the origin lies
in $\st{\circ}{D}_1$. (See Figure 3.1.)
\begin{center}
\unitlength 1.00mm \linethickness{0.4pt}
\begin{picture}(60.00,40.00)
\put(60.00,10.00){\line(-1,0){30.00}}
\put(30.00,10.00){\line(0,1){30.00}}
\put(30.00,40.00){\line(1,-1){30.00}}
\put(30.00,10.00){\circle*{1.33}}
\put(30.00,40.00){\vector(0,-1){18.00}}
\put(34.33,10.00){\line(1,0){25.67}}
\put(60.00,10.00){\line(-1,1){27.00}}
\put(38.00,17.33){\vector(2,1){0.5}}
\bezier{96}(32.33,25.00)(29.00,12.00)(38.00,17.33)
\put(39.67,14.33){\vector(3,1){1.0}}
\bezier{60}(38.00,7.33)(32.67,12.33)(39.67,14.33)
\put(30.00,10.00){\vector(3,2){9.00}}
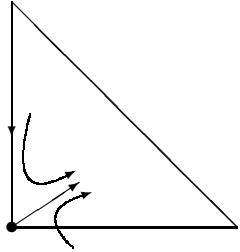
\put(30.00,0.00){\makebox(0,0)[cc]{{\small Figure 3.1. Local
behavior of planar system near the origin when $R_0>1$}}}
\end{picture}
\end{center}
This helps us to find a piece-wise smooth Jordan curve $C$ on
which our vector field $X$ is  either tangent or inward. (See
Figure 3.2.) The Poincar\'e index of such a Jordan curve is 1.
(See \cite{R}, Lemma 5.1.)
\begin{center}
\unitlength 1.00mm \linethickness{0.4pt}
\begin{picture}(60.00,40.00)
\put(60.00,10.00){\line(-1,0){30.00}}
\put(30.00,10.00){\line(0,1){30.00}}
\put(30.00,40.00){\line(1,-1){30.00}}
\put(30.00,10.00){\circle*{1.33}}
\put(50.00,10.00){\vector(-1,1){4.00}}
\put(45.00,25.00){\vector(-1,-1){4.00}}
\put(30.00,40.00){\vector(0,-1){18.00}}
\put(34.33,10.00){\line(1,0){25.67}}
\put(60.00,10.00){\line(-1,1){27.00}}
\bezier{132}(32.67,37.00)(30.67,11.67)(35.33,18.00)
\put(35.50,18.33){\line(2,-3){3.00}}
\bezier{44}(38.67,13.67)(33.33,12.00)(38.67,10.00)
\put(30.00,10.00){\vector(4,3){7.00}}
\put(40.00,3.00){\makebox(0,0)[cc]{{\small Figure 3.2. The Jordan
curve $C$ }}}
\end{picture}\\
\end{center}
We choose this Jordan curve so that it contains all rest points in
$\st{\circ}{D}_1$. If there are two rest points in
$\st{\circ}{D}_1$, they are obtain by a transverse intersection of
the line (3-3) and the conic section $r'=0$ and hence both are
nondegenerate. Thus their Poincar\'e index must be $\pm 1$ which
contradicts $I_X(C)=1$. Therefore there is a unique rest point in
$\st{\circ}{D}_1$ which attracts $D_1-\{i=0\}$ by Theorem 2.1. It
remains to prove that this rest point is hyperbolic. Suppose the
contratry, then it must be nondegenerate by Corollary 3.4. Thus
it is obtained by a tangent (non-transverse) intersection of the
line (3-3) and the conic section $r'=0$. With a slight
perturbation in $\gamma$, we will have either two transverse
intersection in $\st{\circ}{D}_1$ or nothing. ($\gamma$ appears
only in the coefficient of $r$ in $r'=0$ and does not appear in
(3-3) and $R_0$). This is a contradiction with the uniqueness of
the rest point in $\st{\circ}{D}_1$ proven above. $\square$

Theorem 3.1. and Theorem 3.5. provide a complete global analysis
of the planar system (3-1),(3-2) in $\st{\circ}{D}_1$. Since the
dynamics of this system in $\st{\circ}{D}_1$ is equivalent to the
dynamics of proportions system $(2-1)'-(2-3)'$ in $D$, we have
proved the following result which gives a complete global
analysis of the proportions system in the feasibility region $D$.

\pa{Theorem 3.6.} Consider the   proportions system
$(2-1)'-(2-3)'$.
\begin{itemize}
 \item[(i)] If $R_0\leq 1$, then the disease free equilibrium
proportions $(1,0,0)$ is globally asymptotically stable in $D$.
\item[(ii)] If $R_0>1$, then there is a unique rest point
$(s^*,i^*,r^*)$ which is globally asymptotically stable in
$D-\{i=0\}$.
\end{itemize}

\section{Analysis of the model equations}
Consider the original
 model equation (2-1)-(2-3) and recall that the population
equation is\break ${N}'=b_0S+(b_1+\B-\var_1)I+(b_2 -\var_2)R-dN.$
Thus
$$\frac{{N'}}{N}=b_0s+(b_1+\B-\var_1)
i+(b_2-\var_2)r-d. \eqno{(4-1)}$$ If $R_0\leq 1$, then $(s,i,r)\lo
(1,0,0)$ by Theorem 3.6., hence $\dfrac{{N'}}{N}\lo b_0-d$.
Moreover if $R_0>1$ and $I>0$, then $(s,i,r)\lo (s^*, i^*,r^*)$,
i.e. the unique rest point in $\st{\circ}{D}_1$, and
$$\frac{{N'}}{N}\lo b_0
s^*+(b_1+\B-\var)i^*+(b_2-\var _2)r^*-d. $$ We define the second
threshold parameter which governs the total population as follows.
$$R_1=\left\{ \begin{array}{ll}
\dfrac{b_0}{d} & \quad \mbox{if} \ R_0\leq 1,\\
\dfrac{b_0 s^*+(b_1+\B)i^*+b_2 r^*}{d+\var_1 i^*+\var_2 r^*} &
\quad \mbox{if} \ R_0>1.
\end{array} \right. $$
Since $I=0$ is invariant with a linear equation, we may assume
that $I>0$. Now from (2-2), we write
$$\frac{{I'}}{I}=(\B_1-d-\var_1-\alpha)+\lambda s.
\eqno{(4-2)}$$ If $R_0\leq 1$, then $\dfrac{{I'}}{I}\lo
\B_1-d-\var_1-\alpha+\lambda$ and if $R_0>1,\ \dfrac{{I'}}{I}\lo
\B_1-d-\var_1-\alpha+\lambda s^*$. So we define the third
threshold parameter which governs the total number of infective
individuals.
$$R_2=\left\{ \begin{array}{ll}\dfrac{\B_1+\lambda}{d+\var_1+\alpha} &
\quad \mbox{if} \
R_0\leq 1,\\
\dfrac{\B_1+\lambda s^*}{d+\var_1+\alpha} & \quad \mbox{if} \
R_0>1.
\end{array} \right. $$
Notice that $d$ does not appear in the proportions  system and
hence $(s^*, i^*, r^*)$ is independent of $d$. The following
results provide a rather complete global  analysis of the model
equations (2-1)-(2-3).

\pa{Lemma 4.1.} If $I(t)\leq M$ for every $t\geq t_0$, then $R(t)\lo 0$ and if
$I(t)\lo \infty$, then $R(t)\lo \infty$.

\pa{Proof.} From (2-3), we have
${R'}(t)=(\B_2+\alpha)I(t)-(d+\var_2+\gamma) R(t)$. If $I(t)\leq
M$ for $t\geq t_0$, then ${R'} (t)\leq
M(\B_2+\alpha)-(d+\var_2+\gamma)R(t)$ and by Granvell's
inequality \cite{P}, $R(t)\leq M(\B_2+\alpha)
e^{-(d+\var_2+\gamma) (t-t_0)}$ which follows that $R(t)\lo 0$.
Now suppose that $I(t)\lo \infty$. Then by (2-3), we have
${R'}(t)+(d+\var_2+\gamma)R(t)=(\B_2+\alpha)I(t)$ which implies
that $\dfrac{d}{dt} (R(t)e^{(d+\var_2+\gamma)t})=(\B_2+\alpha)I(t)
e^{(d+\var_2+\gamma)t}$. For every
 $M\in \mathbb{R}^+$, there is a $t_0\in \mathbb{R}$
such that $I(t)>M$ for  $t>t_0$. Thus
$$R(t)e^{(d+\var_2+\gamma)t}-R(t_0)
e^{(d+\var_2+\gamma)t_0}=\int_{t_0}^t (\B_2+\alpha) I(t)
e^{(d+\var_2+\gamma)t}> M(\B_2+\alpha)\int_{t_0}^{t}
e^{(d+\var_2+\gamma)t} dt.$$  $$\Rightarrow R(t)>
\left(R(t_0)-\frac{M(\B_2+\alpha)}{d+\var_2+\gamma}\right)
e^{-(d+\var_2+\gamma)(t-t_0)}+\frac{M(\B_2+\alpha)}{d+\var_2
+\gamma}.$$ Since $e^{-(d+\var_2+\gamma)(t-t_0)}$ goes to zero as
$t\lo \infty$, $R(t)> \dfrac{M(\B_2+\alpha)}{2(d+\var_2+\gamma)}$
for large values of $t$. It means that $R(t)\lo \infty$. $\square$

\pa{Theorem 4.2.} (i) If $R_1>1$, then $N(t)\lo \infty$ and if $R_1<1$, then
$N(t)\lo 0$.\\
(ii) If $R_2>1$, then $(I(t),R(t))\lo (\infty,\infty)$ and if $R_2<1$, then
$(I(t),R(t))\lo (0,0)$.

\pa{Proof:} First suppose that $R_0\leq 1$ which implies that
 $(s,i,r)\lo (1,0,0)$ by Theorem 3.6.
Thus $\dfrac{{N'}}{N}\lo b-d$ by (4-1) and $\dfrac{{I'}}{I}\lo
\B_1-d-\var_1-\alpha+\lambda$ by (4-2). If $R_1<1$, then $b-d<0$,
hence $N(t)\lo 0$. Similarly if $R_1>1$, then $b-d>0$, hence
$N(t)\lo 0$. Furthermore If $R_2<1$, then
$B_1-d-\var_1-\alpha+\lambda<0$, hence $I(t)\lo 0$ and by the
above lemma $R(t)\lo 0$. Similarly if $R_2>1$, then $\B_1 - d
-\var_1-\alpha+\lambda>0$, hence $I(t)\lo \infty$ and by the
above lemma $R(t)\lo \infty$.

Now suppose that $R_0>1$. We assumed that $I(t)>0$, hence $i>0$
and  $(s,i,r)\lo (s^*,i^*,r^*)$ by Theorem 3.6. Thus
$\dfrac{{N'}}{N}
 \lo b_0
s^*+(b_1+\B-\var_1)i^*+(b_2-\var_2)r^*-d$ which is positive if
$R_1>1$, hence $N(t)\lo \infty$ and negative if $R_1<1$, hence
$N(t)\lo 0$. Moreover $\dfrac{{I'}}{I}\lo
(\B_1-d-\var_1-\alpha)+\lambda s^*$ which is positive if $R_2>1$,
hence $I(t)\lo \infty$ and then $R(t)\lo \infty$ by the above
lemma. Similarly $\B_1-d\var_1-\alpha+\lambda s^*<0$ if $R_2<1$
and then $(I(t), R(t))\lo (0,0)$. $\square$

We summarize our results in the following table which is the same
as Table 1. in \cite{Bv1}.

\begin{center}
\begin{tabular}{|c|c|c|c|c|c|}
\hline
$R_0$ & $R_1$ & $R_2$ & $N\lo$ & $(s,i,r)\lo$ & $(S,I,R)\lo$\\
\hline
$\leq 1$ & $<1$ & $<1^{\alpha}$ & $0$ & $(1,0,0)$ & $(0,0,0)$\\
$> 1$ & $<1$ & $<1^{\alpha}$ & $0$ & $(s^*,i^*,r^*)$ & $(0,0,0)$\\
$\leq 1$ & $>1$ & $<1$ & $\infty$ & $(1,0,0)$ & $(\infty,0,0)$\\
$\leq 1$ & $>1$ & $>1$ & $\infty$ & $(1,0,0)$ & $(\infty,\infty,\infty)$\\
$>1$ & $>1$ & $>1^{\alpha}$ & $\infty$ & $(s^*,i^*,r)$
 & $(\infty,\infty,\infty)$\\
\hline
\end{tabular}

\vspace*{0.5cm}

$^{\alpha}$ Given $R_0$ and $R_1$, this condition is automatically
satisfied.
\end{center}

\noindent {\bf Conclusion Remarks.} Here we will make some
comparison between the result in \cite{Bv1} and the above results.

1. If we set $b_0=b_1=b_2=b$ and $\B_1=\B_2=\B=0$ in the system
(2-1)-(2-3), we obtain the system (2-1)-(2-3) in \cite{Bv1}.
Moreover by these assumptions, we get $R_0 =
\dfrac{\lambda}{b+\var_1+\alpha}$,
$$R_2=\left\{ \begin{array}{ll}
\dfrac{b}{d} & \quad \mbox{if} \ R_0\leq 1,\\
\dfrac{b}{d+\var_1 i^*+\var_2 r^*} & \quad \mbox{if} \ \ R_0> 1,
\end{array}\right.$$
$$R_3=\left\{ \begin{array}{ll}
\dfrac{\lambda}{d+\var_1+\alpha} & \quad \mbox{if} \ R_0\leq 1,\\
\dfrac{\lambda s^*}{d+\var_1 +\alpha} & \quad \mbox{if} \ R_0> 1.
\end{array}\right.$$
which are the same threshold parameters as in \cite{Bv1}.

2. Comparing their threshold parameters with our ones, we see
that the effect of $b_0$, $b_1$ and $b_2$ appears more clearly
instead of $b$. For example when $R_0>1$, our
 $R_1$ is
$\frac{b_0 s^*+(b_1+\B)i^*+b_2 r^*}{ d+\var_1 i^*+\var_2 r^*}$,
but they obtain $R_1=\frac{b}{d+\var_1 i^*+\var_2 r^*}$ in which
the effects of $b_0$, $b_1$ and $b_2$ are hidden in $b$.

3. There are two vertical transmission parameter, $\B_1$ and
$\B_1$, in our model. The effect of $\B_1$ in $R_0$ and $R_2$ is
crucial. Since $\B_1+\B_2=\B$, we can decrease $\B_1$ by
increasing $\B_2$. Therefore removing more infected newborns
causes more safe situation.

\paragraph{Acknowledgments.} The author would like
to thank Institute for Studies in Theoretical Physics and
Mathematics for supporting this research.

\end{document}